\documentclass[10pt]{article}
\usepackage{amssymb,amsmath,amsfonts} 
\usepackage{mathabx}
\usepackage[colorlinks=true]{hyperref}
\usepackage{comment}

\newtheorem{theorem}{Theorem}

\newtheorem{lemma}{Lemma}

\newtheorem{remark}{Remark}

\newcommand{\diam}{{\rm diam}}

\newcommand{\real}{ {\mathbb R}   }
\newcommand{\torus}{ {\mathbb T}   }
\newcommand{\integer}{ {\mathbb Z}   }
\newcommand{\complex}{ {\mathbb C}   }

\renewcommand{\Im}{\, {\rm Im}\,}

\newcommand\beq[1]{ \begin{equation}\label{#1} }
\newcommand{\eeq}{ \end{equation} }
\newcommand{\beqno}{ \[ }
\newcommand{\eeqno}{ \] }
\newcommand\beqa[1]{ \begin{eqnarray} \label{#1}}
\newcommand{\eeqa}{ \end{eqnarray} }
\newcommand{\beqano}{ \begin{eqnarray*} }
\newcommand{\eeqano}{ \end{eqnarray*} }

\newtheorem{definition}{Definition}
\newcommand\dfn[1]{ \begin{definition}\label{#1} \rm}
\newcommand\edfn{ \end{definition} }
\newcommand{\proof}{\par\medskip\noindent{\bf Proof\ }}
\newcommand\equ[1]{{\rm (\ref{#1})}}
\newcommand{\nl}{{\smallskip\noindent}}
\newcommand{\giu}{{\medskip\noindent}}
\newcommand{\Giu}{{\bigskip\noindent}}
\newcommand{\noi}{{\noindent}}
\newcommand{\qed}{\hskip.5truecm
\vrule width 1.7truemm height 3.5truemm depth 0.truemm
\par\Giu}
\newcommand{\qedeq}{\hskip.5truecm
\vrule width 1.7truemm height 3.5truemm depth 0.truemm}

\newcommand{\e}{\varepsilon}

\renewcommand{\a }{\alpha }
\renewcommand{\b }{\beta }

\renewcommand{\d }{\delta }

\newcommand{\g }{\gamma}
\newcommand{\f }{\varphi}

\renewcommand{\l }{\lambda }

\newcommand{\m }{\mu }

\renewcommand{\t }{\tau }
\renewcommand{\o }{\omega }
\renewcommand{\O }{\Omega }

\def\R{\mathbb R}
\def\T{\mathbb T}

\def\const{{\, c\, }}
\def\dst{\displaystyle}
\def\bks{\, \backslash\, }
\def\meas{{\rm\, meas\, }}
\newcommand\eqby[1]{\stackrel{\equ{#1}}{=}}
\newcommand\leby[1]{\stackrel{\equ{#1}}{\le}}

\newcommand\rball{{\rm B}}
\newcommand\cball{{\,\mathbb B}}
\newcommand\ttM{{\,\mathtt M}}
\newcommand\ttL{{\,\mathtt L}}
\newcommand\ttH{{\,\mathtt H}}
\newcommand\ttD{{\,\mathtt D}}
\newcommand\tth{{\,\mathtt h}}
\newcommand\ttp{\,{\mathtt p}}

\newcommand{\id}{{\mathtt {id}}}

\title{
Explicit estimates on the measure of  primary KAM tori\footnote{
This research was partially supported by the European Research Council under FP7
``Hamiltonian PDEs and small divisor problems: a dynamical systems approach''
and  by MIUR (Italy): PRIN--``Variational and perturbative aspects of nonlinear differential problems''.}
\\
\ \\
\small{(Draft)}
}

\begin{document}

\author{ 
\footnotesize L. Biasco  \& L. Chierchia
\\ \footnotesize Dipartimento di Matematica e Fisica
\\ \footnotesize Universit\`a degli Studi Roma Tre
\\ \footnotesize Largo San L. Murialdo 1 - 00146 Roma, Italy
\\ {\footnotesize biasco@mat.uniroma3.it, luigi@mat.uniroma3.it}
\\ 
}

\maketitle

%
%
%
\begin{abstract}
\noindent
From KAM Theory it follows that the measure of phase points which do not lie on Diophantine, Lagrangian, ``primary''  tori in a nearly--integrable,
real--analytic  Hamiltonian system is $O(\sqrt{\e})$, if $\e$ is the size of the perturbation.  In this paper we discuss how the constant in front of $\sqrt{\e}$ depends on the unperturbed system and in particular on the phase--space domain. 
\end{abstract}
%

\tableofcontents

\section{Introduction}
According to classical KAM theory, the majority of the Lagrangian, non--resonant invariant tori of a ``general'' completely integrable Hamiltonian system persists under the effect of a small enough perturbations (\cite{K}, \cite{A}, \cite{M}; see, also, \cite[\S 6.3]{AKN} for a review and \cite{KAMstory} for a divulgative exposition). 

\nl
Indeed, in bounded regions of the  phase space $\real^n\times\torus^n$ (action--angle) such tori -- which are also called ``primari'' tori -- form a set of positive Liouville (Lebesgue) measure, whose complement has a measure proportional to $\sqrt{\e}$, if $\e$ measures the size of the perturbing function, 
\cite{Nei} and \cite{poschel1982}.

\nl
The square root behavior, in such measure estimates,  is optimal in the sense that, in general, at simple resonances, for $\e\neq 0$, there appear regions of size proportional to $\sqrt{\e}$ free of  primary  invariant tori as trivially shows the example of the simple pendulum with gravity\footnote{Just look at the phase portrait of the simple pendulum $\frac12 p^2 + \e \cos q$, $(p,q)\in\real\times\torus$, and observe that the region enclosed by the separatrix $\frac12 p^2+\e\cos q=\e$ has measure $4\sqrt2\cdot\sqrt\e$.} $\e$.

\nl
It is therefore natural to look for explicit evaluations of the constant in front of $\sqrt\e$ in the KAM measure estimates of the complement of invariant primary tori. 

\nl
In \cite{Nei} and \cite{poschel1982} such constant, which depends on analytic properties of the integrable limit, is left implicit, and, somewhat surprisingly, to the best of our knowledge, there are no explicit evaluations of it in the vast literature on classical KAM theory. On the other hand, KAM is a constructive technique and  discussions about ``KAM constants'' are clearly relevant, as also testified by the large literature on them; compare, e.g., \cite{CG}, \cite{giorgilli}, \cite{cellettiG}, \cite{CC1}, \cite{dlL}, \cite{morbidelli}, \cite{CC2}, \cite{delshamsG}, \cite{CC3}, \cite{cellettiG2}, \cite{CGL}, \cite{CC4}.  

\nl
We also point out that an explicit dependence upon the domain in the above measure estimate is crucial in 
investigating the more complicate problem of the existence and abundance of {\sl secondary tori}, i.e., of those tori which arise by effect of the perturbation around simple resonances. 
In \cite{AKN} it is conjectured that ``in a generic system with three or more degrees of freedom the measure of the non--torus set has order $\e$'', while  
in \cite{BClin} it is given a sketchy proof\footnote{A complete proof will appear elsewhere.} that {\sl the union of primary and secondary tori leave out (for general mechanical systems) a region of measure $\e|\log \e|^a$ (for a suitable $a>0$)}.  To achieve such result one needs to control simultaneously a large number of regions around simple resonances and to apply KAM measure estimates taking into account  different (local) phase--space domains, including neighborhoods of separatrices: To carry out such strategy simple explicit measure estimates --  such as the above -- are necessary.  

\nl
In this paper we compute explicitly the constant in front of $\sqrt{\e}$ up to a constant depending only on $n$ (the number of degrees of freedom) and $\t>n-1$ (the uniform ``Diophantine exponent''). 

\nl
More precisely, we consider a real--analytic, nearly--integrable Hamiltonian 
$$\ttH:(p,q)\in \ttD\times \torus^n\mapsto \ttH(p,q)= h(p)+f(p,q)\in\real$$
where $\ttD$ is an {\sl arbitrary} bounded domain in $\real^n$ and $\torus^n$ is the standard flat $n$--torus (with periods $2\pi$); $f$ is a small perturbation function and the integrable limit $h$ is Kolmogorov non--degenerate on $\ttD$, i.e., its hessian is invertible on $\ttD$. 

\nl
The main result -- Theorem~\ref{KAM} below -- will be formulated in terms of a few (five) parameters, which we now describe briefly (precise definitions will be given in \S~2 below):

\begin{itemize}

\item[\tiny $\bullet$]  The Hamiltonian $\ttH$ is assumed to be real--analytic on $\ttD\times\torus^n$: therefore there exists $r_0>0$ and $0<s\le 1$ such that $\ttH$ is holomorphic on a complex $r_0$--neighborhood of $\ttD$ and a $s$--complex neighborhood of $\torus^n$.

\item[\tiny $\bullet$] The smallness of the perturbation $f$ will be measured by 
$\dst \epsilon:= \frac{\e}{\ttM r_0^2}$ where: $\e=\|f\|_{r_0,s}$ denotes the sup--norm on the above complex neighborhood of $f$ and  
$\ttM:=\|h_{pp}\|_{r_0}$ the sup norm of the Hessian matrix of $h$.

\item[\tiny $\bullet$] The ``torsion'' associated to $h$ will be measured by $\dst \mu:=\frac{\inf_\ttD |\det h_{pp}|}{\ttM^n}$; note that $0<\mu\le 1$ (compare \equ{bove} below). 

\item[\tiny $\bullet$] Last ``independent parameter'' will be the number $\l:=\ttL \ttM$, where $\ttL$ denotes a suitable uniform Lipschitz constant of the local complex inverse of the ``frequency map'' $p\mapsto \o=h_p(p)$ (compare \equ{ivan4}); indeed one can show that $1\le \l\le 2\cdot n! \, \m^{-1}$ (see \equ{misty} below).

\end{itemize}
\nl
Notice that the parameters  $\epsilon$, $\mu$ and $\l$ are  {\sl dimensionless parameters} (i.e., do not have physical dimensions).

\nl
Then, fixed $\nu:=\t+1>n$,  we will show that there exist a positive  constants $c<1$   depending only on $n$ and  $\nu$ such that if the perturbation is so small that
$$
\epsilon\le c \,\frac{\m^6}{\l^2}\ s^{4\nu}\ ,
$$
then one can construct a family ${\cal T}_\a$ of $\ttH$--invariant primary tori. Such tori live in $\ttD_{r_0}\times\torus^n$ (where $\ttD_{r_0}$ is a real $r_0$--neighborhood of $\ttD$) and the $\ttH$--flow on them is analytically conjugated to the Kronecker flow $x\in\torus^n \mapsto x+ \o t$ for a frequency $\o\in\real^n$ which is $(\a,\t)$--Diophantine\footnote{I.e., $|\o\cdot k|\ge \a/|k|_{{}_1}^\t$ for all $k\in\integer^n\bks\{0\}$; compare \equ{raffo} below.}  with $\t=\nu-1$ and  $\a$  proportional to $\sqrt\epsilon$:
$$
\a:=\frac{\l}{ \hat c\, \m\, s^{3\nu}} \, (\ttM r_0)\,\sqrt\epsilon\ ,
$$
where $\hat c<1$ is a suitable constant  depending only on $n$ and  $\nu$.

\nl
The upshot is, then, the following measure estimate (where ``$\meas$" denotes outer Lebesgue measure):
$$
\meas \big((\ttD\times\torus^n)\bks {\cal T}_\a\big) \le  C\, \sqrt\epsilon
$$
where the constant $C$ is given by
$$
C:=\kappa\ 
\big(\max\big\{ \m^2 r_0\, ,\, \diam\, \ttD\big\}\big)^n \cdot \frac{\l^{n+2}}{\m^3\ s^{3\nu}}\ ,
$$
where $\kappa>0$ is a suitable constant  depending only on $n$ and  $\nu$.

\Giu
{\bf Remarks} (i) 
In fact, we shall prove a stronger statement, which is non trivial even in case of $\ttD$ of measure zero or even finite (compare \equ{gusuppone} below). 

\nl
(ii) Of course, more refined estimates are possible if one adds extra hypotheses on the domain $\ttD$ (e.g., smooth boundary) and it would be interesting to give bounds which take into account geometrical properties of $\ttD$. 

\nl
(iii) We do not compute explicitly the dependence upon $n$ (and $\nu$): Indeed it is well known that in such generality explicit bounds on $c$ tend to be quite ``pessimistic'', however, in concrete example, such as a forced pendulum, the standard map or particular three body problems computer--assisted (rigorous) upper bounds on $\e$ are in excellent agreement with experimental data (see, e.g., \cite{CC1}, \cite{dlL}, \cite{CC2},  \cite{CC4} and references therein).

\section{ Notations and set up}

\giu
Given $r>0$, $p_0$ a point of $\real^n$ or $\complex^n$ and $D$ a subset of $\real^n$ or $\complex^n$, we denote:

\beqano
&\rball_r(p_0):=\{p\in\real^n \big| \ |p-p_0|<r\} \ ,  &(p_0\in\real^n)\ , \\
&\cball_r(p_0):=\{p\in\complex^n \big| \ |p-p_0|<r\} \ , & (p_0\in\complex^n)\ , 
\\
&\rball_r(D):=\bigcup_{p_0\in D} \rball_r(p_0) \ , \phantom{AAAAa} & (D\subseteq \real^n)\ ,
\\
&\cball_r(D):=\bigcup_{p_0\in D} {\mathbb B}_r(p_0) \ , \phantom{AAAAa}  & (D\subseteq \complex^n)\ ,
\eeqano
where in $\real^n$ and $\complex^n$, $|x|=|(x_1,...,x_n)|$ will denote the sup--norm $\max_i|x_i|$.

\nl
For a matrix (or a tensor) $A$, $\|A\|$ denotes the standard operator norm $\sup_{|x|=1}|Ax|$. 

\nl
The standard flat $n$--torus $\real^n/(2\pi \integer^n)$ is denoted by $\torus^n$ and, for $s>0$, $\torus^n_s$ denotes its complex neighbourhood of points $q$ with norm of the imaginary part $|\Im q|=|(\Im q_1,...,\Im q_n)|<s$:
$$
\torus^n_s:=\{y\in \complex^n\big| \ |\Im q|<s\}/(2\pi \integer^n)\ .
$$

\nl
If $D$ is an arbitrary bounded set in $\real^n$ and $h$, respectively, $f$,  a real--analytic function (with values in $\real^m$ or in matrix spaces) with bounded holomorphic extension on $\cball_r(D)$ for some $r>0$, respectively,  on $\cball_r(D)\times \torus^n_s $ for some $r,s>0$, we define its analytic sup--norm as, respectively, 
$$
\|h\|_{D,r}:= \sup_{y\in\cball_r(D)}|h(p)|\ ,
\qquad
\|f\|_{D,r,s}:= \sup_{(p,q)\in\cball_r(D)\times \torus^n_s}|f(p,q)|\ .
$$
 The Lipschitz semi--norm of a function $f:\O\to\R^m,$
will be denoted by
$$
|f|_{\rm Lip,\O}:=\sup_{\o_1,\o_2\in\O,\ \o_1\neq\o_2}
\frac{|f(\o_1)-f(\o_2)|}{|\o_1-\o_2|}\,.
$$ 
If $D$ is an open  set and  $H:D\times \torus^n\to\real$ is a $C^2$ function, $\phi_H^t$ denotes its Hamiltonian flow, namely, $\big(p(t),q(t)\big)=\phi_H^t(p,q)$
solves the standard Hamilton equations\footnote{Equivalently, $\phi_H^t$ denotes the Hamiltonian associated to the standard symplectic form $dp\wedge dq=\sum_{i=1}^n dp_i\wedge dq_i$.}
$$
\left\{
\begin{array}{l} 
 \dot p (t) :=\dst  \frac{dp}{dt}(t)= -\partial_q H(p(t),q(t))\\ \ \\
\dot q(t):= \dst \frac{dq}{dt}(t)= \partial_p H(p(t),q(t))
\end{array}
\right. \ ,
\qquad 
(p(0),q(0))=(p,q)\ .
$$
For example,  if $H(p,q)=h(p)$, then the flow $\phi_h^t$ is linear with frequency $\o:=\partial_ph(p)$, namely, $\phi_h^t(p,q)=(p, q+ \o t)$.

\nl
Given  $\a,\t>0$,  a vector $\o\in \real^n$ is said to be $(\a,\t)$--Diophantine if 
\beq{raffo}
|\o\cdot k|:=\big|\sum_{j=1}^n \o_j k_j\big|\ge \frac{\a}{|k|_{{}_1}^\t}\ ,\qquad\forall\ k\in\integer^n\bks\{0\}\ ,
\eeq
where $|k|_{{}_1}:=\sum|k_j|$ denotes the $1$--norm.
It is well know that, fixed $\t>n-1$, almost all (in the sense of Lebegue measure) $\o\in\real^n$ are $(\a,\t)$--Diophantine for some $\a>0$. Indeed such statement  follows immediately observing that\footnote{``meas'' stands for Lebesgue measure, or, in general, for outer  Lebesgue measure.}
\beqno
\meas \big\{\o\in \rball_R(0)\big|\ 
\o \ {\rm is \ not}\  (\a,\t){\rm-\!Diophantine}\big\} \le \const R^{n-1} \a \ ,
\eeqno
with a constant $\const$ depending only on $n$ and $\t$.

\nl
Finally, given a $2n$--vector $(y,x)$, $\pi_1$ and $\pi_2$ denote, respectively, the projectionst on the first and second $n$ components: 
\beq{proiezioni}
\pi_1(y,x)=y\qquad {\rm and}  \qquad \pi_2(y,x)=x\ .
\eeq

\section{Assumptions}\label{assumtions} 

\giu
Fix $n\ge 2$ and $\t>n-1$. Let 
$\ttD$ be any non--empty, bounded subset of $\real^n$. 
Let
\beqno
\ttH:=h+f
\eeqno
with
$h$ and $f$  real--analytic functions with  holomorphic extensions on, respectively,
$\cball_{r_0}(\ttD)$ and $\cball_{r_0}(\ttD)\times \torus^n_s$ for some $r_0>0$ and $0<s\le 1$, and having finite norms:
\beq{bellini}
\ttM:= \|h_{pp}\|_{\ttD, r_0} \,,\qquad\qquad
\e:=\|f\|_{\ttD,r_0,s} \,.
\eeq
Assume that the frequency map $p\in \ttD \to\o= h_p$ is a local diffeomorphism, namely, assume:
\begin{equation}\label{dupa}
d:=  \inf_\ttD|\det h_{pp}| >0\,.
\end{equation}

\section{The local frequency map}\label{frequenze} 
Under assumption \equ{dupa} the frequency map is a {\sl local}  real--analytic diffeomorphism in the neighbourhood of any point of $\ttD$. 
More precisely, the following lemma holds. Define\footnote{Since any eigenvalue of $h_{pp}$ is bounded in absolute value by $\|h_{pp}\|\le \ttM$,  
$d\le \sup_\ttD |\det h_{pp}|\le \ttM^n$.
}
\beq{bove}
\mu:=\frac{d}{\ttM^n}\le 1\ .
\eeq
\begin{lemma}
\label{fichisecchi}
Let
\beq{coco}
c_0=\frac1{8n\cdot n!^2}
\ ,\qquad \qquad 
\hat c_0=\frac1{4n\cdot n!}\ .
\eeq
and define
\beq{rotto}
r_{\star}:= \hat c_0 \m \, r_0\ ,\qquad \qquad
\rho_{\star}:= c_0 \m^2  \ttM\, r_0\ .
\eeq
Then, for every $p_0\in \ttD$  the  frequency map $p \to\o= h_p$ has a real--analytic inverse map, $\o\to \ttp(\o;p_0)$, defined in a neighborhood of $\o_0:=h_p(p_0)$
\begin{equation}\label{ivan2}
\ttp= h_p^{-1}:\ 
\o\in \cball_{\rho_{\star}}(\o_0)\mapsto \ttp(\o;p_0)\in \cball_{r_{\star}}(p_0)\ ,
\end{equation}
with uniform Lipschitz constant\footnote{\label{stewe}
Notice that   on convex domains the Lipschitz semi--norm of a differentiable function coincides with the sup--norm of its Jacobian.} 
\beq{ivan4}
\ttL :=\sup_{p_0\in \ttD}
|\ttp(\cdot;p_0)|_{{\rm Lip},\cball_{\rho_{\star}}(\o_0)}
= \sup_{p_0\in \ttD} \ \sup_{\cball_{\rho_{\star}}(\o_0)} \|  \ttp_\o(\cdot;p_0)\|
\eeq
satisfying
\beq{1936}
\ttL
\le \frac{\ttM^{-1}}{2n\, \hat c_0 \, \m}\ ,
\eeq
and
\begin{equation}\label{XR1000}
\sup_{p_0\in \ttD} \ \ \sup_{\cball_{\frac34\rho_\star}(\o_0)} \|\ttp_{\o\o}(\cdot;p_0)\|
\leq
\frac{\ttL}{\frac{c_0}4 \m^2 \ttM r_0} \ .
\end{equation}
\end{lemma}
\proof
Writing out the inverse of the matrix $h_{pp}$ (Cramer's rule), by Leibniz formula for the $ji$--minor of $h_{pp}$, one has, uniformly on\footnote{Note that $\sup_{i,j}\sup_D|h_{p_ip_j}|\le \ttM:=\sup_D \sup_i \sum_j |h_{p_ip_j}|$.} $D$:
$$
|(h_{pp}^{-1})_{ij}| \le \, \frac1d \, (n-1)!\, \ttM^{n-1}\ ,
$$
which  implies 
\begin{equation}\label{ivan3}
\sup_{\ttD}\| (h_{pp})^{-1}\|\leq \frac{n!}{ \m} \ttM^{-1}\ .
\end{equation}
Let $T:=h^{-1}_{pp}(p_0)$. Then, by standard Cauchy estimates\footnote{\sl If $f:\cball_{r}(D)\to\complex^m$ is holomorphic,
$\partial^\a$ is a partial derivative of order $k=\a_1+\cdots+\a_n$ and $0<r'<r$, then 
$$\sup_{\cball_{r'}(D)}|\partial^\a f|\le \frac{\sup_{\cball_{r}(D)}|f|}{(r-r')^k}\ .
$$},  it follows that for any $p\in \complex^n$ such that $|p-p_0|\le r_{\star}$ one has
\beqno
\|I-T h_{pp}(p)\|\le \|T\|\, \|h_{pp}(p)-h_{pp}(p_0)\|\le \|T\|\, n\, \frac{\ttM}{r_0-r_{\star}}\ r_{\star}\leby{ivan3} \frac{n\cdot n!}{ \m} \frac{r_{\star}}{r_0-r_{\star}}\leby{coco} \frac12\ .
\eeqno
Thus,  by the standard Inverse Function Theorem  (see Appendix~A, Eqn's \equ{IFT1}, \equ{IFT2}, \equ{IFT3})  and Cauchy estimates, relations \equ{ivan2}, \equ{ivan4}, \equ{1936} 
and \equ{XR1000} follow immediately with the constants in \equ{coco}. 
\qed

\nl
For later use, we point out that\footnote{Indeed:  $1=\|I\|=\|h_{pp}(p) h_{pp}^{-1}(p)\|=\| h_{pp}(p) \ttp_\o(h_p(p))\|\le \| h_{pp}(p)\|\ \| \ttp_\o(h_p(p))\| \le \ttM \ttL$.}
\beq{verrecchia}
\l:=\ttL \ttM \ge 1\ ,
\eeq
and that, by \equ{1936},
\beq{misty}
\l\le \frac1{2n\hat c_0}\, \frac1\m= 2\cdot  n!\, \frac1\m\ .
\eeq

%

\section{The classical analytic KAM Theorem}

\begin{theorem}\label{KAM}
Let the assumptions in {\rm \S~\ref{assumtions}} hold and let $\l$, $\mu$ and $c_0$ be as in {\rm \S~\ref{frequenze}} and let $\nu=\t+1$. \\
There exist positive constants $c_\star<1/(8\cdot n!)$ and $\kappa$, depending only on $n$ and $\t$, such that, if  
\beq{baggio}
c:= \frac{c_\star^2}{2^{17}\cdot n^2\, (n!)^6}\ ,\qquad \qquad
\hat c:= \frac{c_\star}{2^4}\ ,
\eeq
and if $\e$ is such that  
\beq{enza} 
\epsilon:= \frac{\e}{\ttM r_0^2}\le c \,\frac{\m^6}{\l^2}\ s^{4\nu}\ , 
\eeq
then the following holds. 
Define
\begin{equation}\label{nicaragua}
\a:=\frac{\l}{ \hat c\, \m\, s^{3\nu}} \, (\ttM r_0)\,\sqrt\epsilon\ ,\qquad
\hat r:= \frac{c_0}2\, \m^2 r_0\ ,\qquad r_\epsilon:=\frac{\l }{c_\star}\, \sqrt\epsilon\, r_0\ .
\end{equation}
Then, there exists a positive measure 
set ${\cal T}_\a\subseteq \rball_{2\hat r}(\ttD)\times \torus^n$ formed by ``primary''  Kolmogorov's tori; more precisely, 
for any point $(p,q)\in{\cal T}_\a$, $\phi^t_\ttH(p,q)$ covers densely an $\ttH$--invariant, analytic, Lagrangian torus,  with $\ttH$--flow analytically conjugated to a linear flow with $(\a,\t)$--Diophantine frequencies 
$\o=h_p(p_0)$, for a suitable $p_0\in\ttD$; each of such tori is a graph over $\torus^n$ $r_\epsilon$--close
to the unperturbed  trivial graph  $\{(p,\theta)=(p_0,\theta)|\ \theta\in \torus^n\}$.\\
Finally, 
the Lebesgue outer measure of $(\ttD\times\torus^n)\bks {\cal T}_\a$ is bounded by:
\beq{gusuppo}
\meas \big((\ttD\times\torus^n)\bks {\cal T}_\a\big) \le  C\, \sqrt\epsilon
\eeq
with
\beqno
C:=\kappa\ 
\big(\max\big\{ \m^2 r_0\, ,\, \diam\, \ttD\big\}\big)^n \cdot \frac{\l^{n+2}}{\m^3\ s^{3\nu}}\ ;
\eeqno
indeed, there exist $N$ and, for $1\le i\le N$,  $p_i\in \ttD$, such that $\ttD\subseteq  \bigcup_{i=1}^N \rball_{\hat r}(p_i)$ and
\beq{gusuppone}
\meas \Big( \big( \bigcup_{i=1}^N \rball_{\hat r}(p_i) \times\T^n\big)\bks  {\cal T}_\a\Big)\le C\,  \sqrt\epsilon\ .
\eeq
\end{theorem}

\Giu
{\bf 6  Remarks}

\begin{itemize}
\item[(i)] The  constant $c_\star$ is, essentially, the ``smallness'' constant appearing in a local KAM normal form (see Theorem~\ref{capra} below). 
The constant $\kappa$ is  given in \equ{cappone}. 

\item[(ii)] Estimate \equ{gusuppone} implies at once \equ{gusuppo} and notice that \equ{gusuppone} is meaningful also in the case of sets $\ttD$ of measure zero (such as a singleton). 

%
\end{itemize}

\section{Proof of KAM Theorem~\ref{KAM}}
The proof of Theorem~\ref{KAM} is divided in six steps.

\subsection{Local reduction}
The first step consists in covering $\ttD$ with $N$ balls centered at points of $\ttD$ (with an explicit upper bound on $N$), thus reducing the Theorem to the special case in which the domain is a ball.   Indeed, the following simple result holds.

\begin{lemma}{\bf (Covering Lemma)} \label{coperta}
Let $E\subseteq\real^n$  be a non--empty set of finite diameter. Then, for any $r>0$ there exists an integer $N$, with\footnote{$[x]$ denotes the integer--part (or ``floor'') function $\max\{n\in \integer |\, \ n\le x\}$, while
$\lceil x\rceil$ denote the ``ceiling function''
$\min\{n\in \integer |\, \ n\ge x\}$.
}
\beq{pippa}
1\le N \le \Big(\Big[ \frac{\diam\, E  }{r}\Big]+1\Big)^n\ ,
\eeq
and $N$ points $p_i\in E$ such that
\beq{zanzara}
E\subseteq \bigcup_{i=1}^N \rball_r(p_i)\ .
\eeq
\end{lemma}
\proof
Let $\d:=\diam\, E$ and
let $z_i=\inf\{x_i|\ x\in E\}$. Then $E\subseteq K:= z+[0,\delta]^n$. Let $0<r'<r$ close enough to $r$ so that $\lceil\d/r'\rceil=[\d/r]+1=:M$. Then, one can 
cover $K$ with $M^n$ closed, contiguous cubes $K_j$, $1\le j\le M^n$, with edge of length $r'$. Let $j_i$ be the indices such that $K_{j_i}\cap E\neq \emptyset$ and pick a $p_i\in K_{j_i}\cap E$;
let $1\le N\le M^n$ be the number of such cubes.
Observe that, since we have chosen the sup--norm in $\real^n$, one has $K_{j_i}\subseteq \rball_{r}(p_i)$ and, therefore,  
 \equ{zanzara} follows with $N$ as in \equ{pippa}. \qed

\nl
We now apply the Lemma with $E=\ttD$ and $r=\hat r$ defined by\footnote{Recall  \equ{rotto},  \equ{bove} and \equ{coco}.}
\begin{equation}\label{dracula}
\hat r :=\frac{\rho_{\star}}{2\ttM}=
\frac{c_0}{2}\, \m^2r_0
\le  \frac{r_0}{128}<
\frac{r_0}{2}\,.
\end{equation}
Thus, Lemma~\ref{coperta} yields that:

\nl
{\sl For suitable $N$ points $p_i\in \ttD$, one has} 
\beq{rufus}
\ttD\subseteq \bigcup_{i=1}^N \rball_{\hat r}(p_i)\ ,\qquad 1\le   N \le \Big(\Big[ \frac{\diam \ttD  }{\frac{c_0}2\m^2 r_0}\Big]+1\Big)^n\ .
\eeq
Notice that, by \equ{dracula}, $\ttH$ is holomorphic and bounded on 
$\cball_{r_0/2} \big(\rball_{\hat r}(p_i)\big)\times\T_s^n$, for every $i\le N$.

\nl
Next we shall prove a ``local'' version of Theorem \ref{KAM}

\subsection{A KAM local normal form apr\`es \cite{poschel2001}}
{\sl Fix one of the balls $\rball_{\hat r}(p_i)$ in the covering \equ{rufus}. We first prove Theorem~\ref{KAM} with $\ttD$ and $r_0$ replaced, respectively, by} 
$$
D_i:=\rball_{\hat r}(p_i)\qquad {\rm and}\qquad \frac{r_0}2\ . 
$$
We shall use  a ``KAM normal form with parameters''; more specifically, we shall use  Theorem~B of \cite{poschel2001}, whose statement we recall here for convenience of the reader.

\nl
Let $r,\a,\tth>0$, $0<s\le 1$, $\t>n-1$ and let $\O\subset \real^n$ be a bounded open set with piecewise smooth boundary; let 
\beqno
\O_\a:=\big\{
\o\in\O, 
\ \ {\rm s.t.}
\ \ {\rm dist}(\o,\partial \O)\geq \a \ \ {\rm and}\ \ 
 \o \ {\rm is}\ (\a,\t)-\!{\rm Diophantine}
\big\}\,;
\eeqno
let $\o\to e(\o)$ and $(I,\theta,\o)\to P(I,\theta,\o)$ be real--analytic functions with holomorphic extension on, respectively, 
$\cball_\tth(\O_\a)$ and $\cball_r(0)\times\torus^n_s\times \cball_\tth(\O_\a)$. Finally, if  $a>0$, we define the ``action rescaling map'':
\beq{scalo}
R_a(I,\theta):=(I/a, \theta)\ .
\eeq
Consider the Hamiltonian function, parametrized by $\o$,  
$$
H(I,\theta,\o):=N(I,\o)+P(I,\theta,\o)\,, \qquad
\text{where}\ \ \ N:=e(\o)+\o\cdot I\,,
$$
with respect to the standard symplectic form $dI\wedge d\theta$; in particular,  the integrable flow $\phi_N^t$ is given by $\phi_N^t(I,\theta)=(I, \theta + \o t)$.

\nl
In \cite{poschel2001} the following result is proven.

\begin{theorem}\label{capra}
Under the above definitions and assumptions, there exist  constants 
\beq{ronaldigno}
0<c_{\star}< \frac{\hat c_{\star}}{4\cdot n!}< \frac1{8\cdot n!}\ ,
\eeq
depending only on $n$ and $\t$,
such that if 
\beq{montone}
|P|_{r,s,\tth}:=\sup_{\cball_r(0)\times\torus^n_s\times \cball_\tth(\O_\a)} |P|\le c_{\star} \a r s^\nu\ ,\qquad \a s^\nu\le \tth \ , \qquad (\nu:=\t+1)\ ,
\eeq
then, there exist a Lipschitz homeomorphism $\f: \O \righttoleftarrow$ and a family of torus embeddings
$$
\Phi:\torus^n\times\O\to \rball_r(0)\times \torus^n\subseteq \real^n\times\torus^n
$$
such the following holds. 
For every $\o\in \O_\a$, $\Phi(\torus^n,\o)$ is an  invariant torus for $H|_{\f(\o)}:=H(I,\theta , \f(\o))$ and
$$
(\phi_{H|_{\f(\o)}}^t \circ \Phi)(\theta,\o)=\Phi(\theta+\o t,\o)\ .
$$
Moreover, for each $\o\in \O$, $\theta\to \Phi(\theta,\o)$ is real--analytic on $\torus^n_{s/2}$  and if  
$$(\theta,\o)\to \Phi_0(\theta,\o):=(0,\theta)$$ 
denotes the trivial torus embedding, one has,
uniformly on, respectively, $\torus^n_{s/2}\times \O$ and $\O$, the following estimates:
\beqa{pesciolone}
&&|R_{r}(\Phi-\Phi_0)|\,,\ \ \a s^\nu |R_{r}(\Phi-\Phi_0)|_{{\rm Lip}, \O}
\leq \frac{|P|}{\hat c_{\star}\a r s^\nu}\,,
\\
\label{pisciolone}
&&|\f -\id|\,,\ \ \a s^\nu |\f -\id|_{{\rm Lip}, \O}
\leq \frac{|P|}{\hat c_{\star} r}
\,.
\eeqa
\end{theorem}

\newpage

\Giu
{ \bf Remarks} \nopagebreak
\begin{itemize}
\item[(i)] The constants $c_{\star}$ and $\hat c_{\star}$ can be taken equal to, respectively,   $\gamma/2$ and $1/c$ where $\gamma$ and $c$ are the constants appearing in Theorem~A of \cite{poschel2001}, where are not explicitly evaluated. 

\nl
In \cite{BB} an infinite dimensional KAM Theorem\footnote{More precisely Theorem 5.1, case (H3) on page 755
and Remark 5.3 on page 758.} (implying Theorem \ref{capra}) is proved, substituting \eqref{montone}
with the stronger (for $s$  small) condition
$|P|_{r,s,\tth}\le \kappa \a r $, where $\kappa=\kappa(s,n,\tau)
:=\kappa_*^{-c \kappa_*}$ with $\kappa_*=(n+\tau)\ln((n+\tau)/s)$ and where $c>0$ is an absolute constant.

\nl
The numerical relations between $c_\star$ and $\hat c_\star$ 
in \eqref{ronaldigno} are assumed for later convenience (and, obviously, are compatible with \cite{poschel2001}).
\item[(ii)] For simplicity -- and because it would play no r\^ole -- in  \equ{pesciolone} we reported slightly weaker estimates with respect to those appearing in  Theorem~A, where in place of $R_{r}$
there appears the rescaling $W(I,\theta):=(I/r,\theta/s)$ (which means that the estimates on the angle components in \cite{poschel2001} are better by a factor $s<1$ than those in \equ{pesciolone}).
\item[(iii)] Actually, the above Theorem is a synthesis of Theorem A and Theorem B in \cite{poschel2001}. In particular, the final measure estimate in Theorem~B is not reported since the constant (and its dependence upon $\O$) is left implicit.  
\item[(iv)] We point out that from the estimates \equ{pesciolone} {\sl does not follows that $\Phi$ is Lipschitz close to the trivial embedding
$\Phi_0$}. Indeed, if  
$\pi_2$ denotes the projection over the $\theta$--component, taking into account that $\a=O(\sqrt{\epsilon})$ (compare \equ{nicaragua}) , 
from \equ{montone} and \equ{pesciolone} it follows that 
$$|\pi_2 \circ \Phi-\pi_2\circ \Phi_0|_{\rm Lip}\le \frac{c_{\star}}{\hat c_{\star}}\cdot \frac1{\a s^\nu} = O\Big(\frac {1}{\sqrt\epsilon}\Big)\ .$$ 
To overcome this fact  one needs suitable  asymmetric rescalings of action and angle variables.  
\item[(v)] As well known (\cite{poschel1982}, \cite{CG}),  maps constructed via KAM methods are smooth in the sense of Whitney, or, what is the same, have $C^\infty$ extensions: Indeed,  $\f$ and $\Phi$ are $C^\infty(\O)$. 
\end{itemize}

\subsection{Applying the KAM normal form to {\large $\ttH|_{\rball_{\hat r}(p_i)\times \torus^n}$}}\label{caprone}
We now apply Theorem~\ref{capra} to $\ttH$ restricted to  
$$D_i\times \torus^n:=\rball_{\hat r}(p_i) \times \torus^n$$ where $\hat r$ is defined in \equ{dracula} and $p_i$ is one of the points introduced in Lemma~\ref{coperta}. Recall that, by \equ{dracula}, $\hat r<r_0/2$ so that $\ttH$ has holomorphic extension to 
$\cball_{r_0/2}(D_i)\times\torus^n_s$.

\nl
Let 
\beq{calcutta0}
\O^{(i)}:=h_p(D_i)=h_p\big(\rball_{\hat r}(p_i)\big)\ ,\qquad\qquad \tth:=\frac{\rho_{\star}}{4}=\frac{\ttM \hat r}2\ ,
\eeq
and notice that, by \equ{dracula}, 
\beqno
\O^{(i)}\subseteq \rball_{\ttM \hat r}(h_p(p_i))\subseteq \rball_{\rho_{\star}/2}(h_p(p_i))= \rball_{2\!\tth}(h_p(p_i)) \quad
\implies \quad \cball_\tth(\O^{(i)})\subseteq \cball_{\frac34\rho_{\star}}(h_p(p_i))\ ,
\eeqno
which, by Lemma~\ref{fichisecchi}, shows that  $h_p$ has an inverse\footnote{\label{iii}Recall that $\ttp=\ttp(\cdot;p_i)$ and therefore depends upon $i$; however for ease of notation we do not indicate explicitly the dependence upon $i$.} $\ttp=h_p^{-1}$ with holomorphic extension
\begin{equation}\label{calcutta}
\ttp=\ttp(\cdot;p_i): \quad\cball_\tth(\O^{(i)}) \to \cball_{r_{\star}}(p_i)
\,.
\end{equation}
Following \cite{poschel2001}, we introduce $\o\in \cball_\tth(\O^{(i)})$ as parameter,  
and let\footnote{In general,  $I$ is complex. Notice that by \equ{enza}, \equ{verrecchia}, \equ{bove} (and the assumption $s\le 1$), $\epsilon< 1/16$ since $c<1/16$.}
\begin{equation}\label{handel}
p=\ttp(\o)+I\ ,\qquad {\rm with}\qquad 
\o\in\cball_\tth(\O^{(i)})\ ,\qquad |I|<r:= \sqrt{\frac{\e}{\ttM}}= \sqrt{\epsilon} \cdot r_0 < \frac{r_0}4\ ,
 \end{equation}
and define
\beq{5101957}
\left\{
\begin{array}{ll} 
N(I,\o):= e(\o)+\o\cdot I:= h(\ttp(\o))+\o\cdot I\\ \ \\ 
\dst P(I,\theta,\o):=\int_0^1 (1-t) h_{pp}(\ttp(\o)+t I) I\cdot I dt\ +f\big(\ttp(\o)+I,\theta)\ ,
\end{array}
\right.
\eeq
so that
$$\ttH(\ttp(\o)+I,\theta)= h(\ttp(\o)+I)+f(\ttp(\o)+I,\theta)=N(I,\o)+P(I,\theta,\o)=:H(I,\theta,\o)\ .
$$
By \equ{calcutta}, \equ{handel}, \equ{rotto}, \equ{coco}, one has that if $\o\in \cball_\tth(\O^{(i)})$ and $I\in \cball_r(0)$, then
\beq{dori}
|\ttp(\o)+I-p_i|\le |\ttp(\o)-p_i|+| I |<r_{\star}+r=\hat c_0\m r_0+\sqrt{\epsilon} r_0<\frac{r_0}{4}+\frac{r_0}{4}=\frac{r_0}2\ ,
\eeq
so that $\ttp(\o)+I \in \cball_{r_0/2}(p_i)$. Thus, 
by \equ{bellini}, \equ{handel} and \equ{5101957}, 
$H$ is real--analytic with holomorphic extension to 
$\cball_{r}(0)\times \torus^n_s\times \cball_\tth(\O^{(i)})$ with
\begin{equation}\label{fregola}
|P|_{r,s,\tth}\leq 2\e\,.
\end{equation}
Thus, if $\a$ is as in \equ{nicaragua} and $r$ as in \equ{handel},
then,  
\beq{fornelletto}
\frac{\ |P|_{r,s,\tth}}{c_\star \a r s^\nu}\le \frac{2\e}{c_\star \a r s^\nu}=\frac{2\hat c}{c_\star}\, \frac{\m s^{2\nu}}{\l} 
\eqby{baggio}\frac18\, \frac{\m s^{2\nu}}{\l}
\le  \frac18\ ,
\eeq
and the first condition in \equ{montone} is satisfied. 
Observe that, 
\beqno
\a s^\nu= \frac{\l}{\hat c \m s^{2\nu}}\, (\ttM r_0) \sqrt\epsilon\leby{enza}\frac{\sqrt{c}}{\hat c} \m^2 \ttM r_0\ ,
\eeqno
and, 
if $\tth$ is as in \equ{calcutta0},  one has, in view of \equ{baggio}, 
\beqno
\frac{\a s^\nu}{\tth}
\leby{enza} \frac{4 \sqrt{c}}{\hat c\, c_0} \le 1\ .
\eeqno
Thus, also the second condition in \equ{montone} is satisfied and
we can apply Theorem~\ref{capra}, obtaining the family of torus embedding\footnote{\label{natravolta}Obviously, also $\Phi$ depends on $i$ but, as above 
(compare footnote \ref{iii}), for ease of notation we do not indicate explicitly the dependence upon $i$.} 
$$\Phi:\torus^n\times\O^{(i)}\to \rball_r(0)\times \torus^n\,\qquad (r=\sqrt\epsilon\, r_0)$$ 
as described in Theorem~\ref{capra}.

\subsection{Kolmogorv's tori: The sets ${\cal T}_\a^{(i)}$ and ${\cal T}_\a$}\label{colmi}
The tori we obtained in the preceding section  live in the ``local'' phase space $\{(I,\theta)|\ (I,\theta)\in \rball_r(0)\times\torus^n\}$. To translate the invariant tori into  the original phase space $\{(p,q)| \ (p,q)\in \rball_{r_0/2}(D_i)\times \torus^n\}$,  we define the $\o$--family of torus embeddings\footnote{Recall \equ{dori} and footnote \ref{natravolta}.}
\beq{beddo}
\theta\mapsto\Psi(\theta,\o)
:=\big(\ttp \circ\f(\o),0\big)+\Phi(\theta,\o)\ , \qquad \o\in\O^{(i)}\ ,
\eeq
which, as function of $\theta$,  is real--analytic on $\T^n_{s/2}$. Then, from \S~\ref{caprone} it follows that  
{\sl for $\o\in\O^{(i)}_\a$ the torus $\Psi(\torus^n,\o)$ is invariant for the flow of $\ttH$ and, furthermore}:
\beqno
(\phi_\ttH^t \circ \Psi)(\theta,\o)=\Psi(\theta+\o t,\o)\ .
\eeqno
We therefore obtain the following family of ``Kolmogorov's tori'' (recall, \equ{dracula},  that $\hat r<r_0/2$):
\beqno
{\cal T}_\a^{(i)}:=\Psi(\torus^n\times \O^{(i)}_\a)\subseteq \rball_{r_0/2}(D_i)\times \torus^n\ ,\qquad\quad
{\cal T}_\a:= \bigcup_{i=1}^N {\cal T}_\a^{(i)}\subseteq \rball_{r_0}(\ttD)\times \torus^n\ ,
\eeqno
Below, we shall show that actually ${\cal T}_\a$ lives in a smaller neighborhood of $\ttD$.

\nl
Analytic quantitative properties of the torus embedding $\Psi$, and hence of the family of Kolmogorov's tori, will be described in detailed in the following section.

\subsection{Properties of the torus embedding}

\begin{lemma}\label{schiaccianoci} Let $\Psi$ be defined as in \equ{beddo} and let $\Psi_0$ denote the   ``trivial embedding''
\beq{bombay}
\Psi_0 \, :\, (\theta,\o)\in\T^n_s\times \cball_\tth(\O^{(i)}) \ \mapsto
\ (\ttp(\o),\theta)\in \cball_{r_\star}(p_i)\times \T^n_s\,.
\eeq
Then\footnote{Recall, \equ{proiezioni},  that $\pi_1$ denotes projection onto the first $n$ components.},
\beq{bue}
\sup_{\O^{(i)}} \sup_{\T^n}|\pi_1\big(\Psi-\Psi_0\big)|\le r_\epsilon\stackrel{\equ{nicaragua}}{:=}\frac{\l \sqrt{\epsilon}}{c_\star} r_0\, \ ,
\eeq
and, hence,
\beqa{corna1}
&&{\cal T}_\a^{(i)} \subseteq \rball_{r_\epsilon}(D_i)\times \T^n=\rball_{\hat r+ r_\epsilon}(p_i)\times \torus^n\ ,\\
&&
\dst {\cal T}_\a\subseteq\bigcup_{i=1}^N\rball_{\hat r+ r_\epsilon}(p_i)\times \torus^n\subseteq \rball_{2 \hat r}(\ttD)\times \T^n
\ .\label{corna2}
\eeqa
\end{lemma}
\proof
By the definitions \equ{beddo} and \equ{bombay}, one has
\beq{moscardini}
\Psi 
=\Psi_0 +
\big(\ttp\circ \f -\ttp ,0\big)
+\Phi -\Phi_0 \,.
\eeq
Notice that,
since $\tth=\rho_\star/4$, from \equ{ivan4} and \equ{XR1000}, one has
\beq{cello2}
\sup_{\cball_\tth(\O^{(i)})} \|\ttp_\o\|
\leq
\ttL  \,,\qquad \qquad
\sup_{\cball_\tth(\O^{(i)})} \|\ttp_{\o\o}\|
\leq
\frac{4}{c_0}\, \frac{\ttL}{\mu^2 \ttM r_0} \ .
\eeq
Thus,
\equ{cello2}, \equ{pisciolone}, \equ{fregola}, one gets
\beq{cello2.5}
\sup_{\O^{(i)}} | \ttp\circ \f - \ttp| \le 
\ttL\, \frac{2 \e}{\hat c_\star r}
\eqby{handel}\frac{2}{\hat c_\star}\ \l r\ .
\eeq
Then, by \equ{pesciolone}, \equ{fornelletto}, \equ{ronaldigno} and \equ{baggio}, one gets \equ{bue}.

\nl
Since (recall \equ{bombay}) 
$$\Psi_0(\T^n\times \O_\a^{(i)})=\ttp(\O_\a^{(i)})\times \T^n\subseteq D_i\times\T^n=\rball_{\hat r}(p_i)\times \T^n\ ,$$
\equ{corna1} follows from \equ{bue}; \equ{corna2} follows since $r_\epsilon\le \hat r$.
\qed

\nl
To control Lipschitz norm we introduce suitable partial rescalings.
Let
\begin{equation}\label{messico}
1<\b:= \frac{\l}{\check c \m s^\nu}\ ,\qquad {\rm with}\qquad \check c := \frac{n \hat c_0 \hat c_\star}{16}=\frac{\hat c_\star}{2^6\cdot n!}
\ .
\end{equation} 
Define, for any $a>0$,  the following   rescaling 
\beqno
S_a: (\theta,x)\in \torus^n\times \real^n \mapsto (\theta,\o):=  (\theta,a x)\ .
\eeqno
Now, recall \equ{scalo} and let
\beq{chiocciola}
\left\{
\begin{array}{ll}
\tilde\Phi:=R_{\b r} \circ \Phi\circ S_\a \\ \ \\ 
\tilde\Phi_0:=R_{\b r} \circ \Phi_0\circ S_\a
\end{array}
\right.\ ,\qquad \qquad
\left\{
\begin{array}{ll}
\tilde\Psi:=R_{\b r} \circ \Psi\circ S_\a \\ \ \\ 
\tilde\Psi_0:=R_{\b r} \circ \Psi_0\circ S_\a
\end{array}
\right.
\eeq
which are defined on the domain $\torus^n_{s/2}\times \frac{1}{\a} \O^{(i)}$. The rescaled version of  \equ{moscardini}, 
then becomes:
\beq{cornoinglese}
 \tilde \Psi-\tilde \Psi_0
 =
 \Big(
 \frac{1}{\b r}(\ttp \circ\f -\ttp)\circ S_\a,0
 \Big)
 +
  \tilde \Phi-\tilde \Phi_0\ .
\eeq
Finally, let
\begin{equation}\label{BWV552}
\Psi_*:=(\tilde\Psi-\tilde\Psi_0)
\circ\tilde\Psi_0^{-1}\, ,
\end{equation}
which is defined on $\frac1{\b r}{D_i}\times\T^n_{s/2}$.
%
%

\nl
The above rescaled embeddings may, now,  be shown to be
close, {\sl in  Lipschitz norm}, to the unperturbed rescaled embeddings:

\begin{lemma}\label{zinedine}
The following bounds hold:
\beqa{cello}
&\dst
\sup_{\torus^n_{s/2}\times \frac{1}{\a} \O^{(i)}}|\tilde \Phi-\tilde \Phi_0|\le \frac{s^\nu}{8}\ ,
&\qquad
\sup_{\torus^n_{s/2}}| \tilde \Phi-\tilde \Phi_0|_{{\rm Lip},  \frac{1}{\a} \O^{(i)}}\le
\frac{s^\nu}{8}\ ,
\\
&\dst
\sup_{\torus^n_{s/2}\times \frac{1}{\a} \O^{(i)}}|\tilde\Psi -\tilde\Psi_0|\leq
\frac{s^\nu}{4}\ ,
&\qquad
\sup_{\torus^n_{s/2}}|\tilde\Psi -\tilde\Psi_0|_{{\rm Lip},  \frac{1}{\a} \O^{(i)}}
\leq \frac14\ ,
\label{oboe}\\
&\dst
\sup_{\frac1{\b r}{D_i}\times\T^n_{s/2}}|\Psi_*|
\leq \frac{s^\nu}4
\, ,\ \  \ \quad
&\qquad
\sup_{\torus^n_{s/2}}|\Psi_*|_{{\rm Lip},\frac1{\b r}{D_i}}
\leq 
\frac14
\ ,
\label{clarinetto}
\\
&\dst 
\  \sup_{\frac1{\b r}{D_i}\times \T^n}\|\partial_\theta\Psi_*\|
\le
\frac{s^{\nu-1}}{2}\ . \ &\ 
\label{terremoto}
\eeqa
\end{lemma}

\proof
Since $\b>1$,  by \eqref{pesciolone} and \equ{fornelletto} we have  that
\beq{cello1bis}
\sup_{\torus^n_{s/2}\times \frac{1}{\a} \O^{(i)}}|\tilde \Phi-\tilde \Phi_0|\le
\sup_{\torus^n_{s/2}\times \frac{1}{\a} \O^{(i)}}|R_\b^{-1} (\tilde \Phi-\tilde \Phi_0)|=
\sup_{\torus^n_{s/2}\times  \O^{(i)}}|R_r(  \Phi-  \Phi_0)|
\leq 
\frac{2 \hat c}{\hat c_\star}\,  \frac{\mu s^{2\nu}}{\lambda}\le \frac{s^\nu}{8}\ ,
\eeq
last inequality holding  because of the definition of $\hat c$ in \equ{baggio}.

\nl
Analogously, (by \eqref{pesciolone}, \equ{fornelletto} and the definition of $\hat c$) we have  that\footnote{\label{cippalippa} If $f$ is a Lipschitz map defined on $\O^{(i)}$, then $f\circ S_a$ is Lipchitz on $\frac1a\O^{(i)}$
and $|f\circ S_a|_{{\rm Lip}, \frac1a\O^{(i)}}=a|f|_{{\rm Lip},\O^{(i)}}$.
}  
\beqno
\sup_{\torus^n_{s/2}}| \tilde \Phi-\tilde \Phi_0|_{{\rm Lip},  \frac{1}{\a} \O^{(i)}}\le
\sup_{\torus^n_{s/2}}|R_\b^{-1} (\tilde \Phi-\tilde \Phi_0)|_{{\rm Lip},  \frac{1}{\a} \O^{(i)}}
=\a \sup_{\torus^n_{s/2}} |R_r(\Phi-\Phi_0)|_{{\rm Lip},\O^{(i)}}
\leq 
\frac{2 \hat c}{\hat c_\star}\,  \frac{\mu s^{\nu}}{\lambda}
\le \frac{s^\nu}{8}\ ,
\eeqno
which, together with \equ{cello1bis}, proves \equ{cello}.

Now, 
by \equ{cello2.5}, \equ{handel}, \equ{messico}, we get
\beq{trombone}
\sup_{\torus^n_{s/2}\times \frac1\a \O^{(i)}}\Big| \Big(
 \frac{1}{\b r}(\ttp \circ \f -\ttp)\circ S_\a,0\Big)
 \Big|
= \frac1{\b r} \sup_\O^{(i)} | \ttp\circ \f - \ttp|
\le 
\frac{2\l }{\hat c_\star \beta} 
\eqby{messico} \frac{\m s^\nu}{2^5 n!}
< \frac{s^\nu}8\ .
\eeq
The first estimate in \equ{oboe} now follows at once in view of \equ{cornoinglese}, the first inequality in \equ{cello} and \equ{trombone}.

\nl
In order to prove the second estimate in \equ{oboe}, in view of \equ{cornoinglese}, we need to estimate the 
 Lipschitz semi--norm of $(\ttp\circ \f-\ttp)$. 
 Fix $\o,\o'\in\O^{(i)}$ and set $p=\ttp(\o)$ and $p'=\ttp(\o')$. Let also $\g(t):=(p'-p)t+p,$ for 
 $t\in [0,1]$ and\footnote{ The introduction of the lifted curve $\tilde \g\subseteq \O^{(i)}$ to join $\o$ and $\o'$ is due to the fact that, in general, $\O^{(i)}$ is not convex.} ${\tilde \g}:=h_p\circ\g$.
Then\footnote{ By Remark (v) in {\bf 7.2}, $\f$ is differentiable; the differentiability of $\f$ almost everywhere also follows, independently, from Rademacher's Theorem.
Notice also that, if $f$ is a function differentiable (a.e.) on $\O^{(i)}$, then $\sup_\O^{(i)} |\nabla f|\le |f|_{{\rm Lip},\O^{(i)}}$ (a.e.), the equality holding if $\O^{(i)}$ is convex.}, 
 \begin{eqnarray*}
 &&
 \Big|
 \ttp\big(\f(\o')\big)-\ttp(\o')-
 \ttp\big(\f(\o)\big)+\ttp(\o) 
 \Big|  
 \\
&&\qquad
= \ \left|
\left[
\int_0^1\Big(  \ttp_\o \big(
\f\big({\tilde \g}(t)\big)
\big) \f_\o\big({\tilde \g}(t)\big)
-
\ttp_\o \big({\tilde \g}(t)\big)
\Big) h_{pp}(\g(t))
\,dt
\right]
(p'-p)
\right|
\\ 
&&\qquad
\leq\ \ttM |p'-p|
\int_0^1\left|
\ttp_\o \big(
\f\big({\tilde \g}(t)\big)
\big)\f_\o\big({\tilde \g}(t)\big)
-
 \ttp_\o \big({\tilde \g}(t)\big)
\right|
\,dt
 \\
&&\qquad
= 
\ \ttM |p'-p|
\int_0^1\left|
\ttp_\o \big(
\f\big({\tilde \g}(t)\big)
\big)\big(
\partial_\o(\f - \id) \big|_{\tilde \g(t)}
\big)
+\ttp_\o \big(
\f\big({\tilde \g}(t)\big)
\big)
-
 \ttp_\o \big({\tilde \g}(t)\big)
\right|
\,dt
\\
&&\qquad
\leby{cello2} \ 
\l |\o'-\o|
\Big( 
\ttL
|\f-\id|_{\rm Lip}
\, +\,
\frac{4\ttL}{c_0\mu^2 \ttM r_0}
|\f-\id|
\Big)\\
&&\qquad
\stackrel{\equ{pisciolone},\equ{fregola}}{\le}
\l |\o'-\o|
\Big( \ttL \frac{2\e}{\hat c_\star \a r s^\nu}+ \frac{4\ttL}{c_0\mu^2 \ttM r_0}\, \frac{2\e}{\hat c_\star r}\Big)\\
&&\qquad
\stackrel{\equ{enza},\equ{handel}}{=} |\o'-\o| \frac{2\l^2  r}{\hat c_\star}\Big( \frac1{\a s^\nu}+\frac{4}{c_0 \mu^2 \ttM r_0}\Big)\ .
\end{eqnarray*}
Now, observe that, by \equ{nicaragua}, \equ{enza}, \equ{baggio}  and \equ{coco} (which implies that $c<\hat c^2 c_0^2/4$), one has
$$
\frac{4}{c_0 \mu^2 \ttM r_0}< \frac1{\a s^\nu}\ .
$$
Thus,
\beqno
|\ttp\circ \f-\ttp|_{{\rm Lip},\O^{(i)}}<
 \frac{4\l^2   r}{\hat c_\star \a s^\nu}
\,,
\eeqno
and, therefore (recalling footnote \ref{cippalippa}),  
\beqano
\frac1{\b r} |(\ttp\circ \f-\ttp)\circ S_\a|_{{\rm Lip},\frac1\a\O^{(i)}}
&=&\frac{\a}{\b r} |\ttp\circ \f-\ttp|_{{\rm Lip},\O^{(i)}}\nonumber\\
&< & \frac{4 \l^2 }{\b \hat c_\star s^\nu} \nonumber\\
&\leby{misty}& \frac{2\l}{n  \hat c_0  \mu \, \hat c_\star s^\nu}\, \frac1\b \nonumber\\
&\stackrel{\equ{messico}, \equ{baggio}}{=}& \frac18\ ,
\eeqano
which, together with the second estimate in \equ{cello}, in view of \equ{cornoinglese}, yields also the second bound in \equ{oboe}.

\nl
To estimate $\Psi_*$ (defined in \equ{BWV552}), observe that 
\beqno
\tilde\Psi_0^{-1}(y,\theta) =\big(\theta,  \a^{-1} h_p(\b r y)\big)\ ,
\qquad\quad 
\tilde\Psi_0^{-1}: \frac1{\b r} {D_i}\times  \torus^n_{s/2}  \stackrel{\rm \small onto}{\to}   \torus^n_{s/2}  \times \frac1\a \O^{(i)}\ .
\eeqno
Thus, the first estimate in \equ{clarinetto} follows immediately from the first bound in \equ{oboe}.
As for the Lipschitz semi--norm of $\Psi_*$, 
by \equ{messico}, \equ{nicaragua},  \equ{handel}, \equ{baggio} and \equ{ronaldigno} 
we have,
for all $\theta\in \T^n_{s/2}$,
that\footnote{$|h_p(\b r\cdot)|_{{\rm Lip},\frac1{\b r}{D_i}}$ denotes the
Lipschitz norm of the function $y\to h_p(\b r y)$ on the rescaled domain 
$(\b r)^{-1} {D_i}$.}
$$
|\tilde\Psi_0^{-1}|_{{\rm Lip},\frac1{\b r}{D_i}}
=\frac1\a|h_p(\b r\cdot)|_{{\rm Lip},\frac1{\b r}{D_i}} 
=\frac{\b r}{\a} |h_p|_{{\rm Lip},{D_i}}
\leq \frac{\b r \ttM}{\a}= \frac{\hat c s^{2\nu}}{\check c}
=\frac{c_\star}{\hat c_\star} 4\cdot n! s^{2\nu}<  s^{2\nu}\ .
$$
Thus, in view of the second estimate in \equ{oboe}, we have,
for all $\theta\in \T^n_{s/2}$,
\beqno
|\Psi_*|_{{\rm Lip},\frac1{\b r}{D_i}}
\leq 
|\tilde\Psi-\tilde\Psi_0|_{{\rm Lip},\frac1{\a}\O^{(i)}}\,
|\tilde\Psi_0^{-1}|_{{\rm Lip},\frac{1}{\b r}{D_i}}
<
\frac{s^{3\nu}}4\,.
\eeqno
By \equ{clarinetto} and Cauchy estimates we get \equ{terremoto}. \qed

\noi
We shall also need the following

\begin{lemma}
Let\footnote{Recall that $r=\sqrt\epsilon\, r_0$ is defined in \equ{handel}.} 
\beq{pele'}
\rho:= \frac{\b s^\nu}4\, r\ .
\eeq
Then,
\begin{equation}\label{flauto}
\tilde\Psi\circ\tilde\Psi_0^{-1}
\Big( \frac{1}{\b r} {D_i}\times\T^n\Big):=
\tilde\Psi\circ\tilde\Psi_0^{-1}
\Big( \frac{1}{\b r} \rball_{\hat r}(p_i)\times\T^n\Big)
\supseteq \frac{1}{\b r}\rball_{\hat r-\rho}(p_i)\times \T^n\,.
\end{equation}
\end{lemma}
\proof
Since\footnote{Notice that it is $\hat r-\rho>0$: Indeed, 
by \equ{messico} and \equ{handel}, we see that
$\rho= {\l \sqrt\epsilon r_0}/(4 \check c \m)$, 
so that (recalling the definition of $\hat r$ in \equ{dracula}) $\rho<\hat r$  is seen to be equivalent to
$\epsilon < 4\, \check c^2\, c_0^2 \,  {\mu^6}/{\l^2} \stackrel{\equ{messico},\equ{coco}}{=} \hat c_\star^2 \mu^6/(2^{16} n^2 n!^6\l^2)$, 
which is guaranteed by \equ{enza}, observing that, by \equ{baggio}, $c<\hat c_\star^2/(2^{16} n^2 n!^6)$.} 

$$\frac{1}{\b r}\rball_{\hat r-\rho}(p_i) = \rball_{\frac{\hat r-\rho}{\b r}}\Big(\frac{p_i}{\b r}\Big)\eqby{pele'}
\rball_{\frac{\hat r}{\b r}- \frac{s^\nu}4}\Big(\frac{p_i}{\b r}\Big)
\ ,
$$
one sees that 
\equ{flauto} will hold if, 
for any given $(y_0,\theta_0)\in   \frac{1}{\b r}\rball_{\hat r-\rho}(p_i) \times \torus^n$, there exists  a point\footnote{As standard, the overline denotes closure  and observe that $y_0+y_1\in \frac{1}{\b r}\rball_{\hat r}(p_i)$.} 
$$(y_1,\theta_1)\in\overline{\rball_{s^\nu/4}(0)}\times \torus^n$$ 
such that 
\beqno
(y_0,\theta_0)=
\tilde\Psi\circ\tilde\Psi_0^{-1}(y_0+y_1,\theta_0+\theta_1)\eqby{BWV552} (y_0+y_1,\theta_0+\theta_1)+  \Psi_*(y_0+y_1,\theta_0+\theta_1)\ .
\eeqno
Such relation is, in turn, equivalent 
to the fixed point equation
\beq{pescespada}
(y_1,\theta_1)=- \Psi_*(y_0+y_1,\theta_0+\theta_1)\ .
\eeq
%
We shall solve \equ{pescespada} in two steps: 
(i), we prove that there exists a unique function $y_*(\theta)$ such that\footnote{Recall, \equ{proiezioni}, that $\pi_i$ denotes projection: $\pi_1(y,\theta)=y$ and $\pi_2(y,\theta)=\theta$.}
\beq{trittico}
y_*(\theta)= - \pi_1\Psi_*(y_0+ y_*(\theta),\theta_0+\theta)\ ,\qquad \forall\  \theta \in\torus^n\ ,
\eeq
and, (ii), we  show that the map
\beq{doncarlos}
\theta\in\torus^n\mapsto \theta + \pi_2 \Psi_*(y_0+ y_*(\theta),\theta_0+\theta)\in\torus^n\ ,
\eeq
is onto,
guaranteeing that there exists a $\theta_1\in\torus^n$ so that $ \theta_1 + \pi_2 \Psi_*(y_0+ y_*(\theta_1),\theta_0+ \theta_1)=0$. These two facts
will show that $y_1:=y_*(\theta_1)$ and $\theta_1$ are solutions of \equ{pescespada}, proving the claim.

\nl
{\sl Proof of (i)}: Let 
$$X:=\{\theta\mapsto y(\theta) \in C(\torus^n,\real^n)|\ \sup_{\torus^n} |y|\le s^\nu/4\}\ ,\qquad  F(y)(\theta):=-\pi_1\Psi_*\big(y_0+ y(\theta),\theta_0+\theta\big)\ .
$$ 
Then, by the first inequality in  \equ{clarinetto}, $F:X\to X$, and  the second inequality in \equ{clarinetto} shows that $F$ is a contraction from $X$ into $X$. Hence, there exists a unique fixed point $y_*\in X$ satisfing \equ{trittico}. Furthermore, since $\Psi_*$ is real--analytic, so is $y_*$ and, in particular, its Jacobian $\partial_\theta y_*$ satisfies the equation
$$
\Big(\id + \pi_1\partial_y\Psi_*\big(y_0+ y(\theta),\theta_0+\theta\big)\Big) \partial_\theta y_*= -\pi_1 \partial_\theta \Psi_*\big(y_0+ y(\theta),\theta_0+\theta\big)\ ,
$$
which, by Neumann series, by the second inequality\footnote{Recall the second remark in footnote \ref{stewe}.} in \equ{clarinetto} and by \equ{terremoto}, yields
\beq{brian}
\sup_{\torus^n} \|\partial_\theta y_*\|\le \frac1{1-\frac14}\ \frac{s^{\nu-1}}{2}=\frac23 \, s^{\nu-1}\ .
\eeq
{\sl Proof of (ii)}:  Observe that from the standard Contraction Lemma it follows easily that\footnote{Indeed, $G\circ\tilde G=\id$ is equivalent to the fixed point equation $\tilde g= - g\circ (\id +\tilde g)$ and if we let $X$ denote $C(\torus^n,\torus^n)$ endowed with the standard metric $d(h_1,h_2):=\sup_{\torus^n} d_{{}_{\torus^n}}\big(h_1(\theta),h_2(\theta)\big)$ (where $d_{{}_{\torus^n}}$ denotes the standard flat metric on $\torus^n$), one sees immediately  that the assumption implies that the map $h\in X\mapsto  - g\circ (\id +h)\in X$ is a contraction from $X$ to $X$, whose unique fixed point yields $\tilde g$. Furthermore, since $g$ is $C^1$, so is $\tilde g$ and the inequality on the Jacobian of $\tilde g$ follows by Neumann series after having differentiated  the identity  $\tilde g= - g\circ (\id +\tilde g)$.}:

\nl
{\sl If  $g$ is a $C^1(\torus^n,\torus^n)$ map such that $\lambda:=\sup_{\torus^n}\| \partial_\theta g\|<1$, then, the map 
$G: \theta\in\torus^n\mapsto \theta + g(\theta)\in \torus^n$  has a unique   inverse  $\tilde G: \theta\in\torus^n\mapsto \theta + \tilde g(\theta)\in \torus^n$
with $\tilde g\in C^1(\torus^n,\torus^n)$ and $\sup_{\torus^n}\| \partial_\theta \tilde g\|\le \lambda/(1-\lambda)$.
}

\nl
Now, recalling \equ{doncarlos}, to finish the proof is enough to check that the  Jacobian of the map
$$\theta \mapsto \pi_2 \Psi_*(y_0+ y_*(\theta),\theta_0+\theta)$$ 
has (operator) norm strictly smaller than one. But, by the second inequality in \equ{clarinetto}, \equ{brian} and \equ{terremoto}, one has, for any $\theta\in\torus^n$, 
\beqano
\big\|\partial_\theta \pi_2 \Psi_*\big(y_0+y_*(\theta),\theta_0+\theta\big)\big\| &\le& 
\big\|\partial_y \Psi_*\big(y_0+y_*(\theta),\theta_0+\theta\big)\big\|
\big\|\partial_\theta y_*(\theta)\big\| 
+ 
\big\|\partial_\theta\Psi_*\big(y_0+y_*(\theta),\theta_0+\theta\big)\big\|
\\
&\le& \frac14\cdot \frac23 \, s^{\nu-1}+ \frac{s^{\nu-1}}2= \frac23 s^{\nu-1}<1\ . \qedeq
\eeqano

\subsection{Measure estimates}
We first provide measure estimates on $(D_i\times\T^n)\bks{\cal T}_\a^{(i)}=(D_i\times\T^n)\bks\Psi(\torus^n\times \O^{(i)}_\a)$. 

\nl
Clearly\footnote{The dot over union denotes ``disjoint union''.},
\begin{equation}\label{arpa}
(D_i\times\T^n)\bks{\cal T}_\a^{(i)}\ \subseteq \ 
\big( (D_i\times \T^n)\setminus \Psi(\T^n\times \O^{(i)}) \big)
\ \dot \cup \ \Psi\big(\T^n \times(\O^{(i)}\setminus\O^{(i)}_\a) \big)\,.
\end{equation}
Now, the following estimates hold.

\begin{lemma}{\bf (Measure estimates)} \label{aveninchi}
Recall \equ{coco}, \equ{baggio} and \equ{messico} and define the following constants:
\beq{cappa}
\kappa_1:=(2\pi)^n\ \frac{n c_0^{n-1}}{2 \check c}
\ ,\qquad\qquad \kappa_2:= \Big(\frac{5\pi}2\Big)^n\ ;
\eeq
\beq{spada}
\kappa_3':= 
2\ n^{\frac{n-1}2} \Big(\sum_{k\neq 0} \frac1{|k|_{{}_1}^{\t+1}}\Big) \frac{ c_0^{n-1}}{\hat c}\ ,\qquad
\kappa_3'':= \frac{2n c_0^{n-1}}{\hat c}\ ,
\qquad
\kappa_3:=\kappa_3'+\kappa_3''\ .
\eeq
Then, one has 
\begin{eqnarray}
&&\meas 
\big( (D_i\times \T^n)\setminus \Psi(\T^n\times \O^{(i)}) \big)
\leq \kappa_1\,  \mu^{2n-3} r_0^n \, \sqrt{\epsilon}
\,,
\label{arpa1}
\\
&&{\rm \meas} 
\big(\Psi(\T^n\times \O^{(i)}\setminus\O^{(i)}_\a) \big)
\leq \kappa_2\, L^n \meas (\O^{(i)}\setminus\O^{(i)}_\a)
\,,
\label{arpa2}
\\
&&\meas (\O^{(i)}\setminus\O^{(i)}_\a)
\leq \kappa_3\,   \frac{\m^{2n-3}\l^2}{s^{3\nu}}\ (\ttM r_0)^n\ \sqrt\epsilon\ .
\label{arpa3}
\end{eqnarray}
\end{lemma}

\proof  
Observe that by \equ{chiocciola}
\beq{cicala}
R_{\b r}\circ \Psi= \tilde\Psi
\circ\tilde\Psi_0^{-1}\circ R_{\b r}\circ \Psi_0\ .
\eeq
Thus, since $\Psi_0(\T^n\times \O^{(i)})=D_i\times \T^n$, we have
\beq{500R}
R_{\b r} \circ \Psi(\T^n\times \O^{(i)})=\tilde \Psi\circ \tilde\Psi^{-1}_0\Big(\frac1{\b r} D_i\times \torus^n\Big)\ .
\eeq
Therefore\footnote{\label{magnifico}Recall that $D_i=\rball_{\hat r}(p_i)$; in the last inequality use that for every $0<x<1$ and for every integer $n>1$, one has  $1-(1-x)^n < nx$ and for the last equality recall \equ{dracula}, \equ{pele'}, \equ{messico}, \equ{handel}.},
\begin{eqnarray*}
{\rm \meas} 
\Big( (D_i\times \T^n)\setminus \Psi(\T^n\times \O^{(i)}) \Big)
&=&(\b r )^n 
\meas \Big( R_{\b r}\big((D_i\times \T^n)\big)\setminus R_{\b r}\circ \Psi(\T^n\times \O^{(i)}) \Big)\\
&\eqby{500R}&(\b r )^n 
{\rm \meas} 
\Big( \big(\frac{1}{\b r}D_i\times \T^n\big)
\setminus \tilde\Psi\circ\tilde\Psi_0^{-1}\big(\frac{1}{\b r}D_i\times \T^n\big) \Big)
\\
&
\stackrel{\eqref{flauto}}\leq &
(\b r )^n 
{\rm \meas} 
\Big( \big(\frac{1}{\b r}D_i\times \T^n\big)
\setminus 
\big(
\frac{1}{\b r}\rball_{\hat r-\rho}(p_i)\times \T^n
\big)
 \Big)\\
&=&
{\rm \meas} 
\Big( \big(D_i\setminus \rball_{\hat r-\rho}(p_i)\big) 
\times \T^n\Big)
\\
&=& (2\pi)^n \Big((2\hat r)^n-\big(2(\hat r -\rho)\big)^n\Big)
\\
&\le &
 (2\pi)^n n 2^n \hat r^{n-1} \rho\\
 &\eqby{cappa}&
\kappa_1 \l \mu^{2n-3} r_0^n \sqrt\epsilon\ ,
\end{eqnarray*}
proving \eqref{arpa1}.  

\Giu
To prove \eqref{arpa2}, observe that if $A\subseteq \O^{(i)}$, from \equ{cicala} and the identity (recall \equ{bombay})
$$
R_{\b r}\circ \Psi_0(\T^n\times A)= \frac1{\b r} \ttp (A)\times 
\T^n\ ,
$$
there follows, by \equ{chiocciola},
\beq{ZADA}
\Psi(\T^n\times A)= R_{\b r}^{-1} \circ \tilde \Psi\circ \tilde \Psi_0^{-1}\Big( \frac1{\b r} \ttp(A)\times \T^n\Big)\ .
\eeq
Observe also that,
since\footnote{Recall the definition of $\Psi_*$  in \eqref{BWV552}.} 
$$\tilde\Psi
\circ\tilde\Psi_0^{-1}
=\id + \Psi_*\ ,$$ 
from \eqref{clarinetto} there follows
\beqno
|\tilde\Psi
\circ\tilde\Psi_0^{-1}|_{{\rm Lip},
\frac1{\b r} D_i\times\T^n}\leq 5/4\ .
\eeqno
Now, for every measurable set $A\subseteq \O^{(i)}$,  one has\footnote{In the first inequality, we use (twice) the following fact: {\sl If $A\subseteq$ is a measurable set and $f:A\to\R^n$ is a Lipschitz map, then $\meas f(A)\le |f|_{{\rm Lip},A}^n \meas A$}.}
\beqano
{\rm meas}\big( \Psi(\T^n\times A) \big)
&\eqby{ZADA}&
(\b r)^n {\rm meas}\Big(
\tilde\Psi
\circ\tilde\Psi_0^{-1}
\big(\frac{1}{\b r}\ttp(A)\times \T^n \big)\Big)\\
&\leq&
(2\pi)^n 
\big(|\tilde\Psi
\circ\tilde\Psi_0^{-1}|_{{\rm Lip},
\frac1{\b r} D_i\times\T^n}\big)^n\ 
\ttL^n\ 
{\rm meas}(A)\\
&\le & \kappa_2 \ttL^n \meas (A)\ ,
\eeqano
and \equ{arpa2} follows. 

\Giu
To prove  \eqref{arpa3},
observe that
\beq{mennea}\O^{(i)}\bks \O^{(i)}_\a\subseteq \{\o\in \O^{(i)}|\ \o \ {\rm is \ not\ (\a,\t)\!-\!Diophantine}\}\cup   \O^{(i)}(\a)
\eeq 
where
$$
\O^{(i)}(\a):=\Big\{
\o\in\O^{(i)}, 
\ \ {\rm s.t.}
\ \ {\rm dist}(\o,\partial \O^{(i)})< \a \ 
\Big\}\,.
$$
Let us begin with estimating the measure of the first set in the r.h.s. of  \equ{mennea} keeping track of constants.
Notice that, if $\hat \O^{(i)}$ denotes the euclidean ball of center $h_p(p_i)$ and radius $\sqrt{n} \ttM \hat r$, then 
$$\O^{(i)}=h_p\big(\rball_{\hat r}(p_i)\big) \subseteq \hat \O^{(i)}\ .
$$
Thus, denoting by $\|\cdot\|$ the euclidean norm in $\R^n$, we have
\beqa{lietome}
\meas \{\o\in \O^{(i)}|\ \o \ {\rm not\ (\a,\t)\!-\!Dioph.}\} &\le &
\meas \Big\{\o\in\hat \O^{(i)}|\ \exists\ k\in\integer^n, \ k\neq 0:\ |\o\cdot k|<\frac{\a }{|k|_{{}_1}^\t}\Big\}\nonumber \\
&\le&
2^n\ \sum_{k\neq 0}(\sqrt{n} \ttM \hat r)^{n-1} \frac{\a}{|k|_{{}_1}^{\t+1}}\nonumber\\
&\stackrel{\equ{dracula}, \equ{nicaragua}}{=}& 
\kappa_3'\ \frac{\ttM^n \l \mu^{(2n-3)}}{s^{3\nu}}\ r_0^n \sqrt{\epsilon}\ .
\eeqa
As for the measure of the second set in \equ{mennea}, we observe that 
either  $\hat r\leq\ttL\a$ or
$\hat r>\ttL\a.$
In the first case we have
$$
{\rm meas} (\O^{(i)}(\a))
\leq  {\rm meas} (\O^{(i)})
\leq \ttM^n {\rm meas} (\rball_{\hat r}(p_i))
= 2^n   \ttM^n 
\hat r^{n}
\leq 2^n  \ttM^n \ttL 
\hat r^{n-1}\a
\,.
$$
In the second case, let
$$
\check r:=\hat r-\ttL\a>0\,.
$$
We claim that 
\beq{camomilla}
h_p \big( \rball_{\hat r}(p_i)\setminus
 \rball_{\check r}(p_i)\big) \supseteq \O^{(i)}(\a)\,.
\eeq
Indeed, by contradiction, assume that there exist  $\o=h_p(p)\in\O^{(i)}(\a),$  $\o_*=h_p(p_*)\in\partial \O^{(i)}$ (namely $|p_*-p_i|=\hat r$) with
$|p-p_i|<\check r,$  and
 $|\o-\o_*|<\a.$
 Then 
 $$
 \ttL\a=\hat r-\check r<|p_* -p|\leq \ttL
 |\o_* -\o|<\ttL\a\,,
 $$
 proving \equ{camomilla}. 
 Thus\footnote{Recall footnote \ref{magnifico}.},
 \begin{eqnarray*}
 {\rm meas} (\O^{(i)}(\a))
 &\leq&
 {\rm meas} \Big(h_p \big( \rball_{\hat r}(p_i)\setminus
 \rball_{\check r}(p_i)\big)  \Big)
 \leq \ttM^n 
  {\rm meas}\big( \rball_{\hat r}(p_i)\setminus
 \rball_{\check r}(p_i)\big)
 =2^n \ttM^n  (\hat r^n-\check r^n)
 \\
&< &(2\ttM)^n \ n\hat r^{n-1} \ttL \a
 \end{eqnarray*}
Thus, in either case, by \equ{dracula} and \equ{nicaragua},  we have
\beq{gaudio}
\meas \O^{(i)}(\a) \le (2\ttM)^n \ n\hat r^{n-1} \ttL \a = \kappa_3''
\frac{\ttM^n \l^2 \mu^{(2n-3)}}{s^{3\nu}}\ r_0^n \sqrt{\epsilon}
\eeq
By \equ{mennea}, \equ{lietome} and \equ{gaudio}, we have
\beqno
\meas (\O^{(i)}\bks \O^{(i)}_\a) \le \kappa_3 \frac{\ttM^n \l^2 \mu^{(2n-3)}}{s^{3\nu}}\ r_0^n \sqrt{\epsilon}\ .
\eeqno
Lemma \ref{aveninchi} is proved. \qed

\nl
From \equ{arpa} and Lemma~\ref{aveninchi} there follows 
\beqno
\meas \Big((D_i\times\T^n)\bks{\cal T}_\a^{(i)}\Big)\le (\kappa_1+\kappa_2\kappa_3) \ \m^{2n} \, r_0^n \frac{\l^{n+2}}{\m^3 s^{3\nu}}\, \sqrt\epsilon\ .
\eeqno
Now, since it is
\beqano
(\ttD\times \T^n)\bks {\cal T}_\a &=& (\ttD\times \T^n)\bks \bigcup_{i=1}^N{\cal T}_\a^{(i)}\\
&\subseteq& \bigcup_{i=1}^N (D_i\times\T^n)\bks \bigcup_{i=1}^N{\cal T}_\a^{(i)}\\
&\subseteq& \bigcup_{i=1}^N (D_i\times\T^n)\bks  {\cal T}_\a^{(i)}\\
&=& \bigcup_{i=1}^N (\rball_{\hat r}(p_i)\times\T^n)\bks  {\cal T}_\a^{(i)}\ .
\eeqano
Now, in view of \equ{rufus}, 
one obtains \equ{gusuppone} with
\beq{cappone}
\kappa:=\frac{2^{2n}}{c_0^n}\  (\kappa_1+\kappa_2\kappa_3) \ ,
\eeq
and   \equ{gusuppo} follows at once. \qed



\appendix

\section{The standard quantitative Inverse Function Theorem}
The following is a standard Inverse Function Theorem in $\complex^n$;
the bar over sets denotes closure.

\Giu
{\bf Proposition} 
{\sl Let $f:\overline{ \cball_r(p_0)}\to \complex^n$ be a holomorphic function with invertible Jacobian $f_p(p_0)$ and with $r$ such that
\beq{IFT1}
\sup_{\overline{\cball_r(p_0)}} \|I- f_p^{-1}(p_0) f_p(p)\|\le \delta<1\ . 
\eeq
Then, there exists a unique holomorphic inverse $g$ of $f$ such that 
\beq{IFT2}
g:\overline{ \cball_\rho(\o_0)}\to \overline{\cball_r(p_0)}\ ,\qquad
{\rm with}\qquad  
\rho:=(1-\delta)\, \frac{ r}{\|f_p^{-1}(p_0)\|}
 \ ,  \qquad \o_0:=f(p_0)\ .
\eeq
Furthermore,
\beq{IFT3}
\dst \sup_{\overline{ \cball_\rho(\o_0)}} \| g_\o\|\le \frac{1}{1-\delta}\, \|f_p^{-1}(p_0)\|\ .
\eeq
If $f$ is real--analytic, so is $g$.
}

\Giu
The elementary proof follows by checking that the map $h\mapsto \Phi(h):=h+f_p^{-1}(p_0)\big(f\circ h -\id\big)$ is a contraction on the space of continuous functions from $\overline{ \cball_\rho(\o_0)}$ in  $\overline{\cball_r(p_0)}$. Then, by the Contraction Lemma,  $g=\lim \Phi^n(p_0)$, which also shows, by Weierstrass Theorem on the uniform limit of holomorphic functions, that $g$ is holomorphic (and real--analytic, if so is $f$). The bound \equ{IFT3} is a general fact following from Neumann series: Indeed, if $A$ and $B$ are $(n\times n)$ matrices and $\|I-AB\|\le\delta<1$, then, by Neumann series, $AB$ is invertible and so are $A$ and $B$, furthermore $\|B^{-1}\|\le \|(AB)^{-1}\|\ \|A\|\le
(1-\delta)^{-1}\|A\|$. 
\qed

\end{document}